\documentclass[10pt]{article}
\usepackage{amsmath,amssymb,amsfonts,amsthm,amscd,graphicx,psfrag,epsfig}
\usepackage{color}
\definecolor{blue}{rgb}{0,0,0.7}
\definecolor{red}{rgb}{0.75, 0, 0}
\usepackage{stmaryrd}
\usepackage[titletoc,toc]{appendix}
\usepackage{hyperref}
\hypersetup{colorlinks, linkcolor=blue, citecolor=cyan}
\usepackage[all]{xy}           
\usepackage{marginnote}  
\usepackage{float}
\usepackage{soul}

\newtheorem{theorem}{Theorem}[section]
\newtheorem{lemma}[theorem]{Lemma}
\newtheorem{proposition}[theorem]{Proposition}
\newtheorem{corollary}[theorem]{Corollary}
\newtheorem{conjecture}[theorem]{Conjecture}
\newtheorem{definition}[theorem]{Definition}

\newcommand{\bpf}{\begin{proof}}
\newcommand{\epf}{\end{proof}}
\newcommand{\bs}{\begin{split}}
\newcommand{\es}{\end{split}}
\newcommand{\be}{\begin{equation}}
\newcommand{\ee}{\end{equation}}
\newcommand{\bt}{\begin{theorem}}
\newcommand{\et}{\end{theorem}}
\newcommand{\bd}{\begin{definition}}
\newcommand{\ed}{\end{definition}}
\newcommand{\bp}{\begin{proposition}}
\newcommand{\ep}{\end{proposition}}
\newcommand{\bl}{\begin{lemma}}
\newcommand{\el}{\end{lemma}}
\newcommand{\bc}{\begin{corollary}}
\newcommand{\ec}{\end{corollary}}
\newcommand{\bcon}{\begin{conjecture}}
\newcommand{\econ}{\end{conjecture}}
\newcommand{\la}{\label}

\newcommand{\Q}{{\mathbb Q}}
\newcommand{\C}{{\mathbb C}}

\newcommand{\hra}{\hookrightarrow}
\newcommand{\lra}{\longrightarrow}

\newcommand{\bg}{\begin{equation}\begin{gathered}}
\newcommand{\eg}{\end{gathered}\end{equation}}

\newcommand{\old}[1]{}

\newcommand{\A}{{\rm A}}
\newcommand{\B}{{\rm B}}

\renewcommand{\P}{\mathbb{P}}

\renewcommand{\H}{{\rm H}}
\renewcommand{\P}{{\rm P}}

\newcommand{\bsp}{\begin{split}}
\newcommand{\esp}{\end{split}}
\newcommand{\epr}{\end{proof}}
\newcommand{\bpr}{\begin{proof}}

\renewcommand{\hra}{{\hookrightarrow}}
\input xy
\xyoption{all}

\begin{document}

\setcounter{section}{1}

\date{}  

\title {Calculating Archimedean Height Pairing via generalized cross-ratio}
\author{Alexander B. Goncharov}

\maketitle

\begin{abstract}

 I answer a question which  S. Bloch asked  during the {\it Regulators V} conference in June 2024 in Pisa.

 \end{abstract}

During the last several years, Spencer Bloch emphasized that Archimedean hight pairings between two homologically trivial cycles of complimentary  arithmetic dimensions, see (\ref{cd1})  for the precise condition, or the closely related classes in the biextension, 
should lead to interesting analogs of the classical cross ratio. 
In this note I calculate such a generalized cross-ratio for  planes in  $\P^n$.

\vskip 2mm
 
Let $\A_1, \A_2$ (respectively  $\B_1, \B_2$) be   $a-$dimensional (respectively $b-$dimensional) planes in $\C\P^n$, where 
\be \la{cd1}
a+b+1=n.
\ee
We assume that  $\A_i\cap \B_j= \emptyset$. 
Let us define a number\footnote{The construction below works for any field $k$.}
$$
r(\A_1, \A_2; \B_1, \B_2) \in \C^\times.
$$

Let $V_{n+1}$ be a complex vector space, and $\C\P^n = \P(V_{n+1}-\{0\})$. Denote by $\widehat \A_i$ and $\widehat \B_j$ the subspaces in $V_{n+1}$ of dimensions $a+1$ and $b+1$ 
 projecting onto the planes $\A_i, \B_j$, $i, j=1,2$. Pick  non-zero volume elements $\alpha_i \in {\rm det}({\widehat \A}_i)$ and $\beta_j \in {\rm det}({\widehat \B}_j)$. 
 Since $
  \A_i \oplus  \B_j =  V_{n+1},
 $ 
  we can define 
$$
(\A_i, \B_j) := \alpha_i \wedge \beta_j \in {\rm det}(V_{n+1}).
$$
\bd \la{Def1} We set 
$$
r(\A_1, \A_2; \B_1, \B_2) := \frac{(\A_1, \B_1)(\A_2, \B_2)}{(\A_1, \B_2)(\A_2, \B_1)}\in \C^\times.$$
\ed
Evidently this number does not depend on the choices of the volume forms $\alpha_i, \beta_j$. Clearly one has 
\be \la{MP}
r(\A_1, \A_2; \B_1, \B_2)\cdot r(\A_2, \A_3; \B_1, \B_2) = r(\A_1, \A_3; \B_1, \B_2).
\ee
$$
r(\A_1, \A_2; \B_1, \B_2) = r(\B_1, \B_2; \A_1, \A_2).
$$

Let me recall  the definition of the Archimedean height pairing for a compact $n-$dimensional K\"ahler manifold $M$ \cite{B}, \cite{Bl}. Let $\Delta \subset M\times M$ be the diagonal. 
Denote by ${\cal H}_\Delta\in {\cal A}^{n,n}(M\times M)$ the harmonic representative of its cohomology class.  There exists a Green current 
$
G_{M\times M}\in {\cal A}^{n-1, n-1}(M)$ such that  
$$
(2\pi i)^{-1}\partial \overline \partial G_{M \times M}= \delta_{\Delta} - {\cal H}_\Delta.
$$ 
\bd Given  homologically trivial 
cycles ${\cal A}, {\cal B}\subset M$ with proper intersection, such that ${\rm dim}{\cal A} + {\rm dim}{\cal B} +1 = {\rm dim}M$
 the Archimedean height pairing $\langle {\cal A}, {\cal B}\rangle$ is defined by
\be
\langle {\cal A}, {\cal B}\rangle:= \int_{{\cal A}\times {\cal B}} G_{M \times M}.
\ee
\ed

Any other current $G_{M \times M}$ solving the equation differs by an exact  current, leaving the integrals  the same. 
Consider the currents 
$$
G_{\cal A}:= \int_{{\cal A}}G_{M \times M}, \ \ \ \ G_{\cal B}:= \int_{{\cal B}}G_{M \times M}.
$$
They  satisfy the differential equations
$$
(2\pi i)^{-1}\partial \overline \partial G_{\cal A} = \delta_{\cal A}, \ \ \ (2\pi i)^{-1}\partial \overline \partial G_{\cal B} = \delta_{\cal B}.
$$
Therefore we have
\be
\langle {\cal A}, {\cal B}\rangle = \int_{{\cal B}} G_{\cal A} = \int_{{\cal A}} G_{\cal B}.
\ee

In particular, we have 
$$
\langle\A_1-\A_2, \B_1-\B_2\rangle := \Bigl(\int_{\A_1(\C)} -  \int_{\A_2(\C)}\Bigr)G_{\B_1-\B_2}. 
$$

\bt \la{Th1} One has 
$$
\langle\A_1-\A_2, \B_1-\B_2\rangle = \log |r(\A_1, \A_2; \B_1, \B_2) |.
$$

\et

\bpr  If  $n=1$, the Green function $G_{\cal A}$ for the divisor ${\cal A}=\{0\}-\{\infty\}$ is $\log|z|$. So its evaluation on the divisor ${\cal B}= \{z\}-\{1\}$ is $\log|z|$. 
Assume that $n>1$. Then using the symmetry we can assume without loss of generality that the planes $\A_i$ are of codimension $\geq 2$. 

Take generic hyperplanes $\H_1, \H_2$ containing the  cycles $\A_1, \A_2$ respectively. Their intersection $\H_1\cap \H_2$ is of codimension $2$. Take a generic subspace ${\rm A}\subset \H_1\cap \H_2$ of 
the same dimension as $A_1$. Let  
$$
\B_i':= \H_1 \cap \B_i; \ \ \ \ \ \ \B_i'':= \H_2 \cap \B_i, \ \ \ \ i=1,2.
$$
 Then we have the following   quadruples of planes in the hyperplanes $\H_1$ and   $\H_2$:

$$
(\A, \A_1; \B'_1, \B_2')\subset \H_1; \ \ \ \ \ \ \ \ \ \ (\A_2, \A; \B''_1, \B''_2)\subset \H_2.
$$
Let us use the notation $\langle *,*\rangle_{\H} $ and $r_\H(...)$  working with planes in the projective space $\H$. 
\vskip 1mm

\bl \la{04}The invariant $r(\A_1, \A; \B_1, \B_2) $ 
remains the same if we  pass to the  hyperplane $H_1 \subset \C\P^n$ and replace $\B_1, \B_2$ by their intersections with $\H_1$, and similarly for $\H_2$: 
\be \la{10}
\begin{split}
&r _{\C\P^n} (\A_1, \A; \B_1, \B_2)= r_{\H_1}(\A_1, \A; \B'_1, \B'_2);\\
&r _{\C\P^n} (\A, \A_2; \B_1, \B_2) = r_{\H_2}(\A, \A_2; \B''_1, \B''_2).\\
\end{split}
\ee
\el

\bpr It is sufficient to show that given a hyperplane $\H$ containing subspaces $\A_1, \A_2$, and setting $\B_i'=\H\cap \B_i$, we get 
$
r _{\C\P^n} (\A_1, \A_2; \B_1, \B_2)= r_{\H}(\A_1, \A_2; \B'_1, \B'_2).
$ Take volume forms $\omega_1, ..., \omega_4$ in $\H$, and  a linear functional $f$ in $\C^{n+1}$ annihilating $\H$. We have the following identity, where each fraction is a number in $\C^\times$, given by 
the ratio of two volume forms: 
 $$
\frac{\omega_1 \wedge f}{\omega_2\wedge f} \cdot \frac{\omega_3 \wedge f}{\omega_4\wedge f} =\frac{\omega_1}{\omega_2} \cdot \frac{\omega_3}{\omega_4}.
 $$
 The claim follows immediately from this. \epr

By the induction on the dimension of the projective space $\H_i$ we have 
\be \la{11}
\begin{split}
&\langle\A_1-\A, \B'_1-\B'_2\rangle_{\H_1} = \log |r_{\H_1}(\A_1, \A; \B'_1, \B'_2) |;   \\
& \langle\A-\A_2, \B''_1-\B''_2\rangle_{\H_2} = \log |r_{\H_2}(\A, \A_2; \B''_1, \B''_2) |.\\
\end{split}
\ee

\bl
\be
\begin{split}
&\langle\A_1-\A, \B'_1-\B'_2\rangle_{\H_1}  = \langle\A_1-\A, \B_1-\B_2\rangle_{\C\P^n}; \\
&\langle\A-\A_2, \B''_1-\B''_2\rangle_{\H_2} =  \langle\A-\A_2, \B_1-\B_2\rangle_{\C\P^n}.\\
\end{split}
\ee
\el

\bpr It is sufficient to check the first identity, and it follows easily from the definitions. An alternative proof follows from Lemma \ref{L2} below. \epr
Substituting this to the left hand sides of identities (\ref{11}),  adding the obtained identities and using   the multiplicativity (\ref{MP}) and Lemma \ref{04}, we get:
\be \la{2}
\begin{split}
&\langle\A_1-\A, \B'_1-\B'_2\rangle_{\H_1} + \langle\A-\A_2, \B''_1-\B''_2\rangle_{\H_2} = \langle\A_1-\A_2, \B_1-\B_2\rangle_{\C\P^n}.\\
\end{split}
\ee

Therefore the left hand side in (\ref{2}) is equal to 
\be
\begin{split}
&\langle\A_1-\A, \B'_1-\B'_2\rangle_{\H_1} + \langle\A-\A_2, \B''_1-\B''_2\rangle_{\H_2} \\
= &\langle\A_1-\A, \B_1-\B_2\rangle_{\C\P^n} + \langle\A-\A_2, \B_1-\B_2\rangle_{\C\P^n} \\
= &\langle\A_1-\A_2, \B_1-\B_2\rangle_{\C\P^n}.   \\\end{split}
\ee
 So by the induction the claim boils down to $\C\P^1$, where it is well known. 
  \epr
  \vskip 2mm
  
 Let us lift  the Archimedean height pairing to an invariant in $\C^\times\otimes \Q$. Since odd dimensional cohomology of ${\Bbb C}\P^n$ vanish, we can avoid biextensions and proceed as follows. Assume that the planes $\A_i, \B_j$ do not intersect. Then the well known  arguments provide an exact sequence\footnote{ Indeed, there are maps 
 \be
 \begin{split}
 &H^{2b+1}(\C\P^n - \A_1\cup \A_2, \B_1\cup \B_2; \Q(b)) \lra H^{2b+1}(\C\P^n - \A_1\cup \A_2; \Q(b))=\Q(0). \\
 &H^{2b+1}(\C\P^n, \B_1\cup \B_2; \Q(b))=\Q(1) \lra H^{2b+1}(\C\P^n - \A_1\cup \A_2, \B_1\cup \B_2; \Q(b)). \\
 \end{split}
 \ee
 Here the second equality in the first row follows from the exact sequence  
 $$
 H^{2b}(\C\P^n, \Q(b))\lra H^{2b}(\A, \Q(b)) = \Q^{\oplus 2}\lra H^{2b+1}(\C\P^n-\A, \Q(b)) \lra 0.
 $$
 The first equality in the second row is similar. }
 $$
 0 \lra \Q(1) \lra H^{2b+1}(\C\P^n - \A_1\cup \A_2, \B_1\cup \B_2; \Q(b))\lra \Q(0) \lra 0.
 $$
Considered, say, as a mixed Hodge structure over $\Q$, this extension defines a class 
 $$
 {\rm cl}(\A_1-\A_2, \B_1-\B_2)\in \C^\times\otimes \Q.
  $$
 
 \bt Under the assumptions above, we have 
$$
{\rm cl}(\A_1-\A_2, \B_1-\B_2) = r(\A_1, \A_2; \B_1, \B_2).
$$
\et

\begin{proof} We use the same notation as in the proof of Theorem \ref{Th1}. 
 The   $n=1$ is the classical case. Let us assume that $n>1$.   
 By the induction on the dimension of the projective space $\H_i$ we have 
\be \la{11a}
\begin{split}
&{\rm cl}(A_1-\A, \B'_1-\B'_2)_{\H_1} =  r_{\H_1}(\A_1, \A; \B'_1, \B'_2) ;   \\
& {\rm cl}(\A-\A_2, \B''_1-\B''_2)_{\H_2} =  r_{\H_2}(\A, \A_2; \B''_1, \B''_2).\\
\end{split}
\ee

\bl \la{L2}
\be
\begin{split}
&{\rm cl}(\A_1-\A, \B'_1-\B'_2)_{\H_1}  = {\rm cl}(\A_1-\A, \B_1-\B_2)_{\C\P^n}; \\
&{\rm cl}(\A-\A_2, \B''_1-\B''_2)_{\H_2} =  {\rm cl}(\A-\A_2, \B_1-\B_2)_{\C\P^n}.\\
\end{split}
\ee
\el

\bpr  Consider the map of the cohomology   induced by the embedding $\H_1\hra \C\P^n$: 
\be
H^{2b+1}(\C\P^n - \A_1\cup \A_2, \B_1\cup \B_2; \Q(b))\lra H^{2b+1}(\H - \A_1\cup \A, \B'_1\cup \B'_2; \Q(b)).
\ee
It  induces an isomorphisms on the sub-object $\Q(1)$ and the quotient object $\Q(0)$.  Therefore it is an isomorphism, and thus preserves the extension class. 
The second  identity is similar.\epr

Substituting this to the left hand sides of identities (\ref{11a}),  adding the obtained identities and using   the multiplicativity (\ref{MP}) and Lemma \ref{04}, we get:
\be \la{2a}
\begin{split}
&{\rm cl}(\A_1-\A, \B'_1-\B'_2)_{\H_1} + {\rm cl}(\A-\A_2, \B''_1-\B''_2)_{\H_2} = {\rm cl}(\A_1-\A_2, \B_1-\B_2)_{\C\P^n}.\\
\end{split}
\ee
 
\epr

Finally, it is well known that the Archimedean height pairing is related to the extension class by $$\log |{\rm cl}(\A_1-\A_2, \B_1-\B_2)| = \langle \A_1-\A_2, \B_1-\B_2\rangle.
$$

\paragraph{Acknowledgment} I am grateful to the referee for useful comments. I am grateful to Spencer Bloch  for discussions during the conference, 
and to the organisers for the invitation to Pisa. This work was supported by the NSF grant DMS-2153059.


\begin{thebibliography}{W}
\bibitem[B]{B} Beilinson A.A.:  {\it Height pairing between algebraic cycles}, K-theory, Arithmetic and Geometry, Yu. I. Manin (Ed.), Lect. Notes in Math. 1289, Springer, 1987. 

\bibitem[Bl]{Bl} Bloch S.: {\it Height pairings},  Journal of Pure and Applied Algebra 34 (1984), 119–145. 

\bibitem[Bl2]{Bl2} Bloch S.: {\it Cycles and biextensions}, Contemporary Mathematics 83 (1989), 19–30.

\end{thebibliography}
 \end{document}